\documentclass{amsart}

\usepackage{amsmath}
\usepackage{amsfonts}
\usepackage{amssymb}
\usepackage{mathtools}
\usepackage{graphicx}
\usepackage{wrapfig}
\usepackage{algorithmic}
\usepackage{todonotes}
\usepackage{xcolor}
\usepackage{array}
\usepackage{extarrows}
\usepackage{wasysym}
\usepackage{dutchcal}
\usepackage{tikz}
\usepackage{pgfplots}
\usetikzlibrary{plotmarks,arrows}
\tikzset{arrow/.style={-latex, shorten >= 1ex, shorten <=1ex}}
\usepackage{pgfplotstable}
\pgfplotsset{compat=1.12}
\usetikzlibrary{external}
\usetikzlibrary{shapes.arrows}
\usepackage[algo2e,ruled, linesnumbered]{algorithm2e}
\usepackage{footmisc}
\usepackage{bbm}
\usepackage{subcaption}

\usepackage[skip=0pt, font=footnotesize, width = \linewidth]{caption}
\numberwithin{equation}{section}

\usepackage{hyperref}
\hypersetup{colorlinks=true,allcolors=blue}

\usepackage[capitalize,nameinlink,noabbrev]{cleveref}
\crefname{section}{section}{sections}
\crefname{subsection}{subsection}{subsections}
\Crefname{section}{Section}{Sections}
\Crefname{subsection}{Subsection}{Subsections}
\crefname{algo}{Algorithm}{Algorithms}
\crefname{table}{Table}{Tables}
\crefformat{equation}{\textup{#2(#1)#3}}
\crefrangeformat{equation}{\textup{#3(#1)#4--#5(#2)#6}}
\crefmultiformat{equation}{\textup{#2(#1)#3}}{ and \textup{#2(#1)#3}}
{, \textup{#2(#1)#3}}{, and \textup{#2(#1)#3}}
\crefrangemultiformat{equation}{\textup{#3(#1)#4--#5(#2)#6}}%
{ and \textup{#3(#1)#4--#5(#2)#6}}{, \textup{#3(#1)#4--#5(#2)#6}}{, and \textup{#3(#1)#4--#5(#2)#6}}
\Crefformat{equation}{#2Equation~\textup{(#1)}#3}
\Crefrangeformat{equation}{Equations~\textup{#3(#1)#4--#5(#2)#6}}
\Crefmultiformat{equation}{Equations~\textup{#2(#1)#3}}{ and \textup{#2(#1)#3}}
{, \textup{#2(#1)#3}}{, and \textup{#2(#1)#3}}
\Crefrangemultiformat{equation}{Equations~\textup{#3(#1)#4--#5(#2)#6}}%
{ and \textup{#3(#1)#4--#5(#2)#6}}{, \textup{#3(#1)#4--#5(#2)#6}}{, and \textup{#3(#1)#4--#5(#2)#6}}
\crefdefaultlabelformat{#2\textup{#1}#3}

\definecolor{pantone234} {RGB}{175, 0, 120} 
\definecolor{blue}{rgb}{0.0, 0.2, 0.8}
\definecolor{green}{rgb}{0.0, 0.5, 0.0}
\definecolor{orange}{rgb}{1.0, 0.4, 0}

\ifpdf
\hypersetup{
	pdftitle={Local training and enrichment based on a residual localization strategy},
	pdfauthor={T. Keil, M. Ohlberger, F. Schindler, J. Schleu{\ss}}
}
\fi

\title[Training and enrichment based on a residual localization strategy]{Local training and enrichment based on a residual localization strategy}

\author{Tim Keil, Mario Ohlberger, Felix Schindler, and Julia Schleu{\ss}}
\address{Institute for Analysis and Numerics, University of M\"unster, Einsteinstra\ss e 62, 48149 M\"unster, Germany, {\tt\{tim.keil,mario.ohlberger,felix.schindler,julia.schleuss\}@uni-muenster.de}}

\date{\today}

\thanks{This work was funded by the BMBF under contract 05M20PMA and by the Deutsche Forschungsgemeinschaft (DFG, German Research Foundation) under Germany's Excellence Strategy EXC 2044-390685587, Mathematics M\"unster: Dynamics~--~Geometry~--~Structure.}



\begin{document}
	
	\begin{abstract}
To efficiently tackle parametrized multi and/or large scale problems, we propose an adaptive localized model order reduction framework combining both local offline training and local online enrichment with localized error control. For the latter, we adapt the residual localization strategy introduced in [Buhr, Engwer, Ohlberger, Rave, SIAM J. Sci. Comput., 2017] which allows to derive a localized a posteriori error estimator that can be employed to adaptively enrich the reduced solution space locally where needed. Numerical experiments demonstrate the potential of the proposed approach.
	\end{abstract}

	\maketitle


\section{Introduction}
In this contribution we are concerned with reduced basis model order reduction for parameterized elliptic partial differential equations, see, e.g., \cite{BCOW:17}.
The objective is to compute efficient and certified approximate solutions of parametrized multi and/or large scale problems. Especially if repeated simulations of complex heterogeneous problems for different parameters are required as, for instance, within optimization or inverse problems, a numerical approximation via standard techniques such as the finite element method often becomes prohibitively expensive. To overcome this issue, we propose an adaptive localized model order reduction framework, which completely avoids any expensive global fine scale computations and only requires local fine scale computations for both the construction of local reduced basis functions and the error estimation. Moreover, the approach is easily parallelizable and thus well-suited to be employed on modern computer architectures.

In the offline phase, we precompute a set of carefully chosen problem-adapted local basis functions that serves as an initial reduced basis. To this end, we exploit the methodology of optimal local approximation spaces, see, e.g., \cite{BabLip11,MaSch22,SchSme22,SchSmeMaa23,SmePat16}, and their efficient approximation using techniques from randomized numerical linear algebra \cite{BuhSme18,HalMarTro11}. Next, in the online phase, we use a localized residual-based a posteriori error estimator to investigate the accuracy of the reduced solution for any given new parameter. 
As the error estimator is localized, we can exploit it to adaptively enrich the reduced solution space locally where necessary, cf. \cite{Albetal12,OhlSch15}. The approach thus guarantees the accuracy of reduced solutions given any possibly insufficient reduced basis. The guiding idea allowing for localized error control is to employ and adapt the residual localization strategy proposed in \cite{Buhetal17} which enables to localize the computation of the dual norm of the residual. We also refer to \cite{BleMalVoh20,CiaVoh17} for similar localization results. Moreover, we note that the combination of offline training and online enrichment allows to flexibly balance and shift workload based on, e.g., the respective application or the employed computer architecture. A related approach combining offline training and online enrichment for component-based parametric model order reduction has been recently proposed in \cite{SmeTad22}.

For an overview on localized model order reduction methods we refer to \cite{Buhetal20}. Besides, adaptive (localized) model order reduction methods have recently been successfully employed for solving PDE-constrained optimization problems in \cite{Keietal21,KeiOhl24,KeiOhlSch23} and parameter identification problems in \cite{KarKeiOhl23}.

The remainder of this paper is organized as follows. In \cref{sec:FOM} we introduce the general problem setting and the discretization scheme used at full order level. Subsequently, we propose the locally adaptive model order reduction framework in \cref{sec:LROM} including both the offline training phase outlined in \cref{subsec:offline} and the online enrichment phase with localized error control presented in \cref{subsec:online}. Finally, we showcase numerical experiments in \cref{sec:numerical_experiments} to demonstrate the potential of the proposed approach and draw conclusions in \cref{sec:conclusion}.

\section{Problem setting and full order model}
\label{sec:FOM}

Let $\Omega \subset \mathbb{R}^d$ denote a bounded domain and let $\mathcal{P} \subset  \mathbb{R}^q$ with $q \in \mathbb{N}$ denote a parameter space. As a prototype of a parametrized multi and/or large scale problem, we consider the following diffusion-reaction problem (in weak formulation): For $\mu \in \mathcal{P}$ find $u_\mu \in H^1_0(\Omega)$\footnote{For simplicity of notation we here assume homogeneous Dirichlet boundary conditions on $\partial\Omega$.} such that
\begin{equation} \label{eq:PDE_weak}
\int_\Omega \kappa_\mu \nabla u_\mu \nabla v + r_\mu u_\mu v = \int_\Omega f_\mu v \qquad \text{for all } \,v  \in H^1_0(\Omega).
\end{equation}
Here, $\kappa_\mu \in L^\infty(\Omega)^{d\times d}$ denotes a diffusion tensor and $r_\mu \in L^\infty(\Omega)$ denotes a reaction rate which are both strictly positive, essentially bounded, and of possibly high contrast and multiscale structure. Moreover, $f_\mu \in L^2(\Omega)$ denotes a source term.

We assume that a non-overlapping decomposition of the computational domain~$\Omega$ is given by a coarse grid $\mathcal{T}_H$ with subdomains $T_j \in \mathcal{T}_H$ for $1 \leq j \leq N_H$ as illustrated in \cref{fig:grid_visualizations}(a). Moreover, each coarse subdomain $T_j$ is further decomposed by a local triangulation $\tau_h(T_j)$ that is assumed to resolve all fine scale features of the parametrized multiscale problem \cref{eq:PDE_weak} in order to guarantee accurate numerical approximations. On each coarse subdomain we consider a standard conforming piecewise linear finite element space $V_h(T_j)\coloneqq S^1(\tau_{h}(T_j))$. We emphasize that the proposed approach is not limited to this specific choice and any scheme suitable for discretizing the considered problem could be employed. Furthermore, we couple local spaces in a non-conforming way and thus define the global solution space as $V_h \coloneqq \bigoplus_{j=1}^{N_H} V_h(T_j)$, which results in functions in $ V_h$ being two-valued on coarse inner grid faces.

We then define the full order model via a symmetric weighted interior penalty discontinuous Galerkin scheme w.r.t. the coarse grid: Find $u_{h,\mu} \in V_h$ such that
\begin{equation} \label{eq:PDE_weak_DG}
a_\text{DG}(u_{h,\mu},v_h;\mu) = l_\text{DG}(v_h;\mu)  \qquad \text{for all } \,v_h  \in V_h,
\end{equation}
where the DG bilinear form $a_\text{DG}$ is given by 
\begin{equation*}
a_\text{DG}(v_h,w_h;\mu) \coloneqq \sum_{T\in\mathcal{T}_H} \int_T \kappa_\mu \nabla v_h \nabla w_h + r_\mu v_h w_h + \sum_{\gamma\in\mathcal{F}(\mathcal{T}_H)}a_\text{DG}^\gamma(v_h,w_h;\mu)
\end{equation*}
and the linear form is given by $l_\text{DG} \coloneqq \sum_{T\in\mathcal{T}_H} \int_T f_\mu v_h$. Here, $\mathcal{F}(\mathcal{T}_H)$ denotes the set of all faces of $\mathcal{T}_H$ and the DG coupling bilinear form $a_\text{DG}^\gamma$ for a face $\gamma$ is given by 
\begin{equation*}
a_\text{DG}^\gamma(v_h,w_h;\mu) \coloneqq \int_\gamma \langle \kappa_\mu \nabla v_{h} \cdot n_{\gamma} \rangle  \lbrack w_{h} \rbrack +  \langle \kappa_\mu \nabla w_{h} \cdot n_{\gamma} \rangle  \lbrack v_{h} \rbrack + \frac{\sigma \{ \kappa_{\mu_*}\}}{h_\gamma} \lbrack v_{h} \rbrack \lbrack w_{h} \rbrack.
\end{equation*}
A unique normal $n_\gamma$ pointing away from the adjacent subdomain $T^-$ is assigned to each face $\gamma\in\mathcal{F}(\mathcal{T}_H)$, where an inner face is given by $\gamma = T^- \cap T^+$ and a boundary face is given by $\gamma = T^- \cap \partial\Omega$. Moreover, the weighted average and jump of a function $v_h \in V_h$ across a coarse grid face $\gamma$ are given by $\langle v_h \rangle \coloneqq w_\gamma^- v_h \vert_{T^-} + w_\gamma^+ v_h \vert_{T^+}$ and $\lbrack v_h \rbrack \coloneqq v_h \vert_{T^-} -v_h \vert_{T^+}$ for an inner face and by $\langle v_h \rangle \coloneqq \lbrack v_h \rbrack \coloneqq v_h $ for a boundary face.
For a chosen reference parameter $\mu_*\in \mathcal{P}$ and associated permeability $\kappa_{\mu_*}$, the weights $w_\gamma^-$ and $w_\gamma^+$ are given by $w_\gamma^- \coloneqq \kappa_{\mu_*}\vert_{T^-} (\kappa_{\mu_*}\vert_{T^+} + \kappa_{\mu_*}\vert_{T^-})^{-1}$ and 
$w_\gamma^+ \coloneqq \kappa_{\mu_*}\vert_{T^+} (\kappa_{\mu_*}\vert_{T^+} + \kappa_{\mu_*}\vert_{T^-})^{-1}$. In addition, $\{ \kappa_{\mu_*}\}$ amounts to half the harmonic average of $\kappa_{\mu_*}$ for an inner face and to $\kappa_{\mu_*}\vert_{T^-}$ for a boundary face.
Finally, the penalty parameter $\sigma$ is chosen carefully to ensure coercivity of $a_\text{DG}$ with respect to $\Vert \cdot \Vert_h$, the norm we equip $V_h$ with which is given by
\begin{equation*}
\Vert v_h \Vert_h^2 \coloneqq  \sum_{T\in\mathcal{T}_H} \Vert \kappa_{\mu_*} \nabla v_h \Vert^2_{L^2(T)} + \Vert r_{\mu_*} v_h \Vert_{L^2(T)}^2
+ \sum_{\gamma\in\mathcal{F}(\mathcal{T}_H)}\frac{\sigma \{ \kappa_{\mu_*}\}}{h_\gamma} \Vert \lbrack v_{h} \rbrack \Vert^2_{L^2(\gamma)}.
\end{equation*}

\section{Localized model order reduction}
\label{sec:LROM}

The key idea of the localized model order reduction approach is to exchange the high-dimensional local spaces $V_h(T_j)$ by carefully chosen problem-adapted localized reduced basis spaces $V_\text{rb}(T_j) \subset V_h(T_j)$ of much smaller dimension. After construction, the reduced local spaces are again coupled in a non-conforming fashion and we thus define the global reduced solution space as $V_\text{rb} \coloneqq \bigoplus_{j=1}^{N_H} V_\text{rb}(T_j)$. The global reduced problem is then given by a Galerkin-projection of \cref{eq:PDE_weak_DG} onto $V_\text{rb}$: For $\mu \in \mathcal{P}$ find $u_{\text{rb},\mu} \in V_\text{rb}$ such that
\begin{equation} \label{eq:PDE_weak_reduced}
a_\text{DG}(u_{\text{rb},\mu},v_\text{rb};\mu) = l_\text{DG}(v_\text{rb};\mu)  \qquad \text{for all } \,v_\text{rb}  \in V_\text{rb}.
\end{equation}
To construct appropriate low-dimensional local approximation spaces, we propose a combination of offline training and online enrichment with localized error control. We highlight that in both steps only local fine scale problems need to be solved and no global solves of problem \cref{eq:PDE_weak_DG} are required, neither for the construction of local reduced basis functions nor for the error estimation. Moreover, as the local computations are independent of each other, they are easily parallelizable.

\subsection{Local offline training}\label{subsec:offline}

We propose to perform a first offline training phase, where a set of problem-adapted local basis functions in each of the coarse subdomains is precomputed for some chosen training parameters. These basis functions are then used as an initial reduced basis in the online phase.

For this purpose, we employ the well-established methodology of constructing (quasi-)optimal local approximation spaces. In the following, we briefly outline the main concepts and refer to, e.g., \cite{BabLip11,BuhSme18,MaSch22,SchSme22,SchSmeMaa23,SmePat16} for further details. 

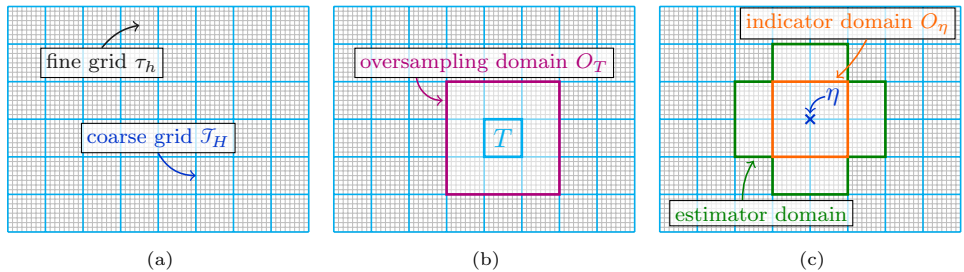
\begin{figure}[t]
	\centering
	\hfill
	\begin{subfigure}{0.315\textwidth}
		\centering
		\begin{tikzpicture}[scale=0.5]
		\draw[step=0.125, black!25] (-1,0) grid (7,6);
		\draw[step=1.0, cyan, line width = 0.6pt] (-1,0) grid (7,6);
		\node[fill=white, fill opacity = 0.75, text opacity = 1, inner sep = 2pt, draw]  at (3,2.5) (TH) { \footnotesize \textcolor{blue}{coarse grid $\mathcal{T}_{H}$}};
		\node[fill=white, fill opacity = 0.75,  text opacity = 1, inner sep = 2pt, draw]  at (1.5,4.5) (Th) {\color{black!90} \footnotesize fine grid  $\tau_{h}$};
		\draw[->, bend right=30, blue, semithick] (3,2.1) to (4,1.5);
		\draw[->, bend left=20, semithick, black!90] (1.5,4.9) to (2.5,5.5);
		\node at (3.05,-0.75) {\scriptsize (a)};
		\end{tikzpicture}
	\end{subfigure}
	\hfill
	\begin{subfigure}{0.315\textwidth}
		\centering
		\begin{tikzpicture}[scale=0.5]
		\draw[step=0.125, black!25] (-1,0) grid (7,6);
		\draw[step=1.0, cyan, line width = 0.6pt] (-1,0) grid (7,6);
		\draw [pantone234, fill=white, fill opacity=0.5,  text opacity = 1, line width = 1pt] (2,4) -- (2,1) -- (5,1) -- (5,4) -- cycle;
		\draw [cyan, line width = 1pt] (3,3) -- (3,2) -- (4,2) -- (4,3) -- cycle;
		\node[text opacity = 1] at (3.5,2.5) {\textcolor{cyan}{${T}$}};
		\node[fill=white, fill opacity = 0.75, inner sep = 2pt, draw, text opacity = 1] at (3,4.5) {\color{pantone234}{\footnotesize oversampling domain ${O_T}$}};
		\draw[->, bend right=40, pantone234, semithick] (1.2,4.1) to (1.9,3.5);
		\node at  (3.05,-0.75) {\scriptsize (b)};
		\end{tikzpicture}
	\end{subfigure}
	\hfill
	\begin{subfigure}{0.315\textwidth}
		\centering
		\begin{tikzpicture}[scale=0.5]
		\draw[step=0.125, black!25] (-1,0) grid (7,6);
		\draw[step=1.0, cyan, line width = 0.6pt] (-1,0) grid (7,6);
		\draw [green, fill=white, fill opacity=0.5, line width = 1pt] 
		(2,4) -- (1,4) -- (1,2) -- (2,2) -- (2,1) -- (4,1) -- (4,2) -- (5,2) -- (5,4) -- (4,4) -- (4,5) -- (2,5) -- cycle;
		\draw [orange, line width = 1pt] (2,4) -- (2,2) -- (4,2) -- (4,4) -- cycle;
		\node at (3.6,3.6) {\textcolor{blue}{${\eta}$}};
		\draw[->, bend right=40, blue, semithick] (3.4,3.6) to (3,3.15);
		\draw [blue, thick] plot [only marks,mark=x, mark options={scale=2.5}] coordinates{(3,3)};
		\node[fill=white, fill opacity = 0.75, draw,  text opacity = 1, inner sep = 2] at (4,5.55) { \footnotesize \textcolor{orange}{indicator domain $O_\eta$}};
		\draw[->, bend left=10, orange, semithick] (4.5,5.12) to (3.7,4.1);
		\node[fill=white, fill opacity = 0.75, inner sep = 2pt, text opacity = 1,  draw] at (1.7,0.5) { \footnotesize \textcolor{green}{estimator domain}};
		\draw[->, bend left=15, green, semithick] (1.2,0.85) to (1.5,1.9);
		\node at  (3.05,-0.75) {\scriptsize (c)};
		\end{tikzpicture}
	\end{subfigure}
	\hfill
	\caption{(a) Non-overlapping domain decomposition into coarse subdomains. (b) Local subdomain~$T$ with associated local oversampling domain $O_T$. (c) Local error indicator domain $O_\eta$ and local estimator domain corresponding to coarse grid node $\eta$.} \label{fig:grid_visualizations}
\end{figure}

\vspace{0.5em}
\textbf{Optimal local approximation spaces.} The key idea and goal is to identify and extract basis functions that are relevant for approximating the local solution space of the PDE. For a given training parameter $\mu_\text{train} \in \mathcal{P}_\text{train}$ in a finite set of training parameters $\mathcal{P}_\text{train}\subset \mathcal{P}$ and a local target subdomain $T\in \mathcal{T}_H$, we thus define a transfer operator $P_{O_T \rightarrow T}^{\mu_\text{train}}$ whose range is the space of local solutions of the PDE \cref{eq:PDE_weak_DG} for $\mu_\text{train}$ on $T$. In more detail, $P_{O_T \rightarrow T}^{\mu_\text{train}}$ 
\begin{enumerate}
	\item takes arbitrary Dirichlet boundary values on the boundary $\partial O_T$ of a larger oversampling domain $O_T$ (typically one additional layer of coarse neighbouring elements, see \cref{fig:grid_visualizations}(b) for an illustration),
	\vspace{2pt}
	\item solves the PDE \cref{eq:PDE_weak_DG} locally on $O_T$ for $\mu_\text{train}$,
	\vspace{2pt}
	\item and restricts the corresponding solution to $T$.
\end{enumerate}
We here refer to, e.g., \cite{BabLip11,MaSch22,SchSme22,SchSmeMaa23,SmePat16}. Since $P_{O_T \rightarrow T}^{\mu_\text{train}}$ is proven to be compact, its singular value decomposition (SVD) can be employed to approximate its range. More precisely, it can be shown that the space $\Lambda_k$ spanned by the $k$ leading left singular vectors of $P_{O_T \rightarrow T}^{\mu_\text{train}}$ is an optimal approximation space in the sense of Kolmogorov, meaning that it minimizes the approximation error among all linear spaces of dimension $k$. Moreover, we have that $\Vert P_{O_T \rightarrow T}^{\mu_\text{train}} - \text{proj}_{\Lambda_k} P_{O_T \rightarrow T}^{\mu_\text{train}}\Vert = \sigma_{k+1}$, where $\text{proj}_{\Lambda_k}$ denotes the orthogonal~projection onto $\Lambda_k$ and $\sigma_{k+1}$ is the $k+1$st singular value\footnote{For simplicity of notation we here omit the dependence of $\Lambda_k$ and $\sigma_{k+1}$ on $T$ and $\mu_\text{train}$.} of $P_{O_T \rightarrow T}^{\mu_\text{train}}$. Due to the fast decay of singular values usually observed in numerical experiments, few left singular vectors are sufficient for an accurate approximation of the local solution space.

\vspace{0.5em}
\textbf{Quasi-optimal local approximation spaces.} To further reduce computational costs and enable efficient parallel computations, we approximate the optimal local space $\Lambda_k$ via random sampling techniques \cite{BuhSme18,HalMarTro11} and thus avoid computing the exact SVD of $P_{O_T \rightarrow T}^{\mu_\text{train}}$. To this end, we apply $P_{O_T \rightarrow T}^{\mu_\text{train}}$ to $k+p$ randomly drawn Dirichlet boundary values, i.e. we solve \cref{eq:PDE_weak_DG} locally on $O_T$ for $\mu_\text{train}$ and random boundary values and restrict the solutions to $T$. Here, $p$ is an oversampling parameter typically not greater than $2$ or $3$. The space $\Lambda_ {k+p}^\text{rand}$ spanned by the $k+p$ resulting local solutions restricted to $T$ yields an approximation that converges provably at a nearly optimal rate of order $\sqrt{k} \,\sigma_{k+1}$. Moreover, based on a probabilistic a posteriori error estimator a so-called \textit{randomized range finder} algorithm has been proposed \cite{BuhSme18} that adaptively constructs a local approximation space $\Lambda_ {m}^\text{rand}$ that satisfies the property $\mathbb{P}(\Vert P_{O_T \rightarrow T}^{\mu_\text{train}} - \text{proj}_{\Lambda_{m}^\text{rand}} P_{O_T \rightarrow T}^{\mu_\text{train}}\Vert \leq \texttt{tol}) > (1 - \varepsilon_{\texttt{fail}})$, where the accuracy \texttt{tol} and the failure probability $\varepsilon_{\tt {fail}}$ are prescribed by the user.\footnote{For simplicity of notation we here omit the dependence of $\Lambda_{m}^\text{rand}$ on $T$ and $\mu_\text{train}$. Moreover, we note that the dimension $m$ of the space $\Lambda_{m}^\text{rand}$ depends on $T$, $\mu_\text{rand}$, \texttt{tol}, and $\varepsilon_{\tt {fail}}$.} For further details we refer to~\cite{BuhSme18}. 

Finally, we define the initial local reduced approximation space $V_\text{rb}(T)$ as the span of basis functions of the spaces $\Lambda_ {m}^\text{rand}$ for all training parameters $\mu_\text{train}$ and set the initial global reduced solution space as $V_\text{rb} \coloneqq \bigoplus_{T\in \mathcal{T}_H} V_\text{rb}(T)$ as introduced above.

\subsection{Local online enrichment with localized error control}\label{subsec:online}

Having constructed an initial reduced basis as outlined above, in the online phase we aim to evaluate the accuracy of reduced solutions for new parameters of interest and enrich the current reduced solution space in case it is not rich enough yet. For this purpose, we use a localized a posteriori error estimator that is derived based on a residual localization strategy \cite{BleMalVoh20,Buhetal17,CiaVoh17} and solve local enrichment problems if necessary, cf. \cite{Albetal12,OhlSch15}. We thus achieve to again only carry out local fine scale computations and completely avoid any computations that scale with the dimension of the global fine grid.

\vspace{0.5em}
\textbf{Localized error control based on residual localization.} In \cite{Buhetal17} a reliable, efficient, and locally computable a posteriori error estimator has been proposed. Its derivation crucially relies on localizing the dual norm of the residual by exploiting orthogonality of the residual with respect to lowest-order shape functions on the coarse grid, more precisely partition of unity functions associated with the coarse grid nodes. While in \cite{Buhetal17} an overlapping decomposition of the computational domain and a conforming coupling of local approximation spaces via a partition of unity approach is considered, we here transfer the result to the non-conforming setting.
For further details we refer to \cite{Buhetal17} and to \cite{BleMalVoh20,CiaVoh17}, where similar results have been shown.

To investigate the accuracy of the reduced solution $u_{\text{rb},\mu_\text{new}} \in V_\text{rb}$ (cf. problem~\cref{eq:PDE_weak_reduced}) for a new parameter of interest $\mu_\text{new} \in \mathcal{P}$, we introduce the residual which is for any parameter $\mu \in \mathcal{P}$ and corresponding reduced solution $u_{\text{rb},\mu}$ given by 
\begin{equation*}
R(u_{\text{rb},\mu};\mu) \in V_h', \quad R(u_{\text{rb},\mu};\mu) \lbrack v_h \rbrack \coloneqq l_\text{DG}(v_h;\mu) - a_\text{DG}(u_{\text{rb},\mu},v_h;\mu). 
\end{equation*} 
To avoid prohibitively expensive computations of the global dual norm of the residual, which are required in classical a posteriori error estimates, we adapt the localization strategy proposed in \cite{Buhetal17}. To this end, we introduce a partition of unity on the coarse grid $\mathcal{T}_H$, where each function $\varphi_\eta^\text{pu}$ is associated to a grid node $\eta \in \mathcal{N}(\mathcal{T}_H)$, and the support of $\varphi_\eta^\text{pu}$ is referred to as indicator domain $O_\eta$, see \cref{fig:grid_visualizations}(c) for an illustration. For the localization result to hold, we moreover need to add $\varphi_\eta^\text{pu}\vert_T$ to the local reduced space $V_\text{rb}(T)$ for all $\eta$ that satisfy $T \subset O_\eta$ for each subdomain $T\in \mathcal{T}_H$. We then obtain the following localized estimator for the approximation error between the solution of the full and the reduced order model (cf. \cref{eq:PDE_weak_DG} and \cref{eq:PDE_weak_reduced}):
\begin{align}
\begin{split}\label{eq:error_estimator}
\Vert u_{h,\mu_\text{new}} - u_{\text{rb},\mu_\text{new}} \Vert_h &\leq \Delta_\text{loc}(u_{\text{rb},\mu_\text{new}};\mu_\text{new}) \\
&\coloneqq \alpha_{\mu_\text{new}}^{-1} \, C_{\text{pu},V_\text{rb}} \Bigg(\sum_{\eta\in\mathcal{N}(\mathcal{T}_H)} \Vert R(u_{\text{rb},\mu_\text{new}};\mu_\text{new}) \Vert^2_{V(O_\eta)'} \Bigg)^\frac{1}{2},
\end{split}
\end{align}
where $\alpha_{\mu_\text{new}}$ denotes the coercivity constant of $a_\text{DG}(\cdot,\cdot\,;\mu_\text{new})$ with respect to $\Vert \cdot \Vert_h$ and $V(O_\eta)$ denotes the restriction of $V_h$ to $O_\eta$.\footnote{A detailed derivation of an upper bound for the partition of unity stability constant $C_{\text{pu},V_\text{rb}} \coloneqq \sup_{v_h \in V_h\backslash \{0\}} \Vert v_h \Vert_h^{-1} \big(\sum_{\eta\in\mathcal{N}(\mathcal{T}_H)} \inf_{v_\text{rb} \in V_\text{rb}\cap V(O_\eta)} \Vert \text{proj}_{V(O_\eta)} (\varphi_\eta^\text{pu} v_h) - v_\text{rb} \Vert_h^2 \big)^{1/2} $ would exceed the scope of this article and is thus subject of a subsequent publication.}

In addition to reliability, we highlight that $\Delta_\text{loc}$ is an efficient error estimator, cf.~\cite{Buhetal17}, and solely locally computable. To compute the dual norm of the residual locally on the indicator domain $O_\eta$, a slightly larger estimator domain as depicted in  \cref{fig:grid_visualizations}(c) is required due to the non-conformity on coarse grid faces and the resulting  coupling terms in the bilinear form $a_\text{DG}$.

As the localized error estimator \cref{eq:error_estimator} initially decomposes into \textit{node-based} local error indicators, we distribute the value of the local dual norms on all supporting coarse subdomains to obtain \textit{element-based} local error indicators:
\begin{equation*}
\delta_{T,\text{loc}}(u_{\text{rb},\mu_\text{new}};\mu_\text{new})^2\mspace{-1mu} \coloneqq\mspace{-1mu} \sum_{\eta\in\mathcal{N}(\mathcal{T}_H), \,T \subset O_\eta}\mspace{-1mu} \frac{1}{ \# \lbrace T \in \mathcal{T}_H \vert T \subset O_\eta\rbrace}  \Vert R(u_{\text{rb},\mu_\text{new}};\mu_\text{new}) \Vert^2_{V(O_\eta)'}.
\end{equation*}

\vspace{0.5em}
\textbf{Local online enrichment.} In case the local error indicator $\delta_{T,\text{loc}}$ derived above indicates that the current reduced solution space $V_\text{rb}$ is not rich enough in a local subdomain $T$, we propose to solve the following local enrichment problem (cf. \cite{Albetal12,OhlSch15}) on the associated oversampling domain $O_T$ (cf. \cref{fig:grid_visualizations}(b)): Find a local correction $\Psi_{T,\mu_\text{new}} \in V_h(O_T)$ such that
\begin{equation*} 
a_\text{DG}\big\vert_{O_T}(u_{\text{rb},\mu_\text{new}} + \Psi_{T,\mu_\text{new}},v_h;\mu_\text{new}) = l_\text{DG}\big\vert_{O_T}(v_h;\mu_\text{new})  \qquad \text{for all } \,v_h  \in V_h(O_T),
\end{equation*}
which amounts to solving the PDE locally on $O_T$ with the residual as right-hand side.
Here, $V_h(O_T)$ denotes the restriction of $V_h$ to $O_T$. Finally, we then enrich the local reduced approximation space $V_\text{rb}(T)$ with the restriction of $\Psi_{T,\mu_\text{new}}$ to $T$. In this way, the local error indicator $\delta_{T,\text{loc}}$ facilitates the design of locally adaptive procedures in which the reduced solution space $V_\text{rb}$ is iteratively enriched locally where needed.

\section{Numerical experiments}
\label{sec:numerical_experiments}

In this section, we demonstrate the potential of the proposed localized model order reduction approach. To this end, we consider a diffusion-reaction problem with highly heterogeneous parametric permeability. The experiments have been performed using \texttt{pyMOR}\footnote{see \url{https://pymor.org}} \cite{MilRavSch16} for the model order reduction as well as \texttt{dune-gdt}\footnote{see \url{https://github.com/dune-gdt/dune-gdt}} and the DUNE framework \cite{Basetal21} for the discretization.

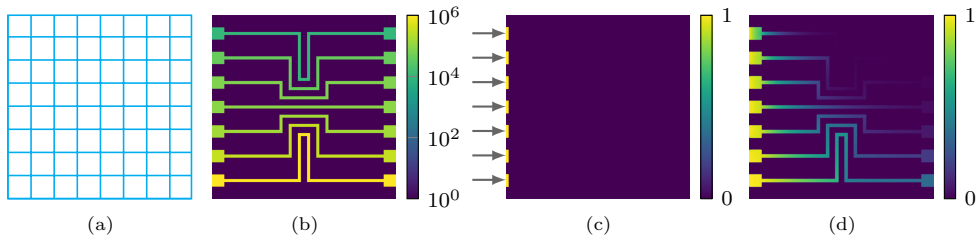
\begin{figure}[t]
	\centering
	\begin{tikzpicture}
	\begin{axis}[
	hide axis,
	width=4.02cm,
	height=4.02cm,
	point meta max=1,
	point meta min=-1,
	xmajorticks=false,
	xmin=0, xmax=8,
	ymajorticks=false,
	ymin=0, ymax=8,
	yshift = -0.2cm
	]
	\end{axis}
	\draw [cyan,semithick] (0,0) rectangle (2.44,2.44);
	\draw[step = 2.46/8, cyan, semithick] (0,0) grid (2.44,2.44);
	\node at (1.23,-0.35) {\scriptsize (a)};
	\end{tikzpicture}
	\hspace{0.05cm}
	\begin{tikzpicture}
	\begin{axis}[
	hide axis,
	colorbar,
	colorbar style={font=\footnotesize,width=0.15cm, ytick={0,2,4,6}, yticklabels={$10^0$,$10^2$,$10^4$,$10^6$}, xshift=-0.15cm},
	colormap name = viridis,
	width=4.02cm,
	height=4.02cm,
	point meta max=6,
	point meta min=0,
	xmajorticks=false,
	xmin=0, xmax=8,
	ymajorticks=false,
	ymin=0, ymax=8,
	]
	\addplot graphics [includegraphics cmd=\pgfimage,xmin=0, xmax=8, ymin=0, ymax=8] {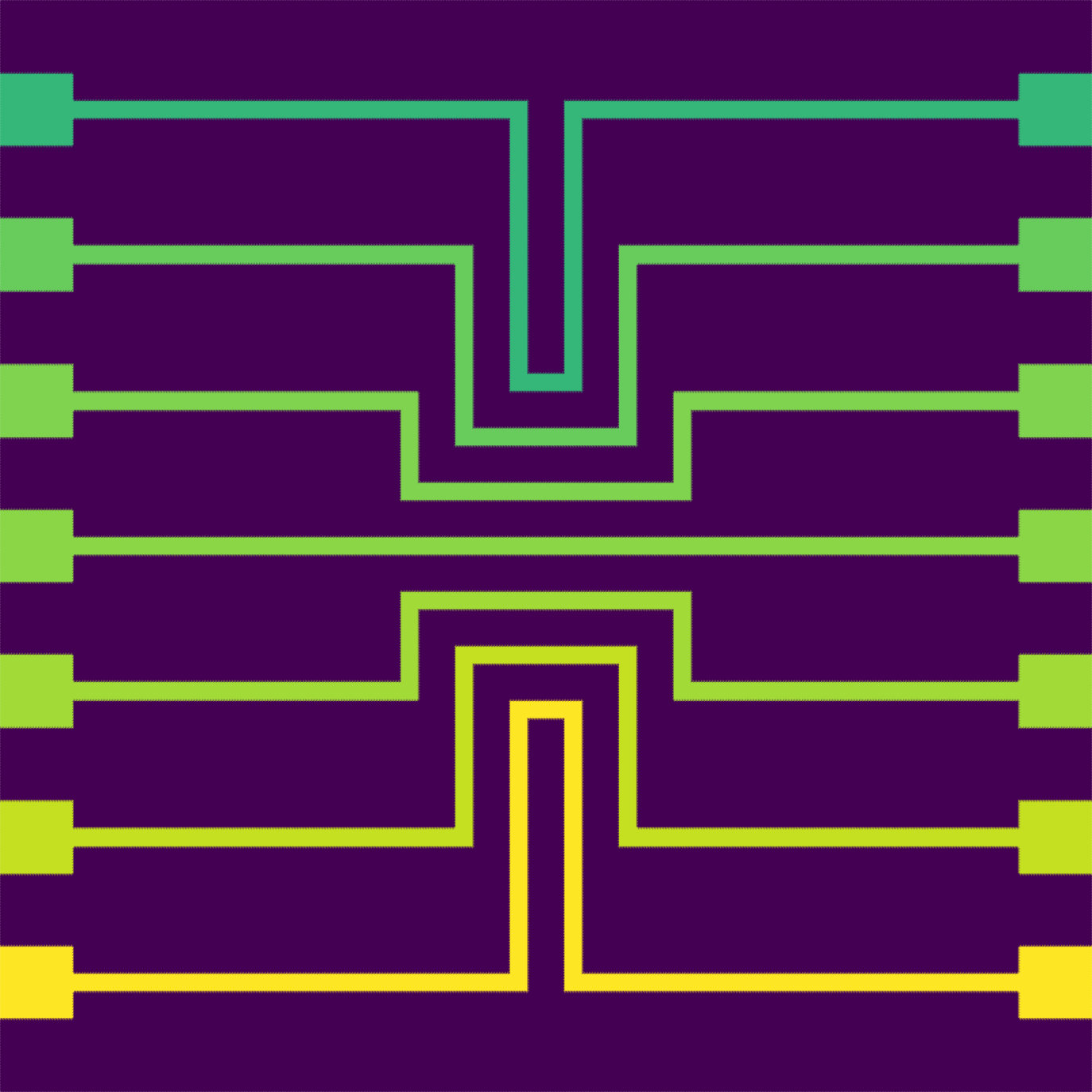};
	\end{axis}
	\node at (1.23,-0.35) {\scriptsize (b)};
	\end{tikzpicture}
	\hspace{-0.35cm}
	\begin{tikzpicture}
	\begin{axis}[
	hide axis,
	colorbar,
	colorbar style={font=\footnotesize,width=0.15cm,ytick = {0,1}, xshift=-0.15cm},
	colormap name = viridis,
	width=4.02cm,
	height=4.02cm,
	point meta max=1,
	point meta min=0,
	xmajorticks=false,
	xmin=0, xmax=8,
	ymajorticks=false,
	ymin=0, ymax=8,
	]
	\addplot graphics [includegraphics cmd=\pgfimage,xmin=0, xmax=8, ymin=0, ymax=8] {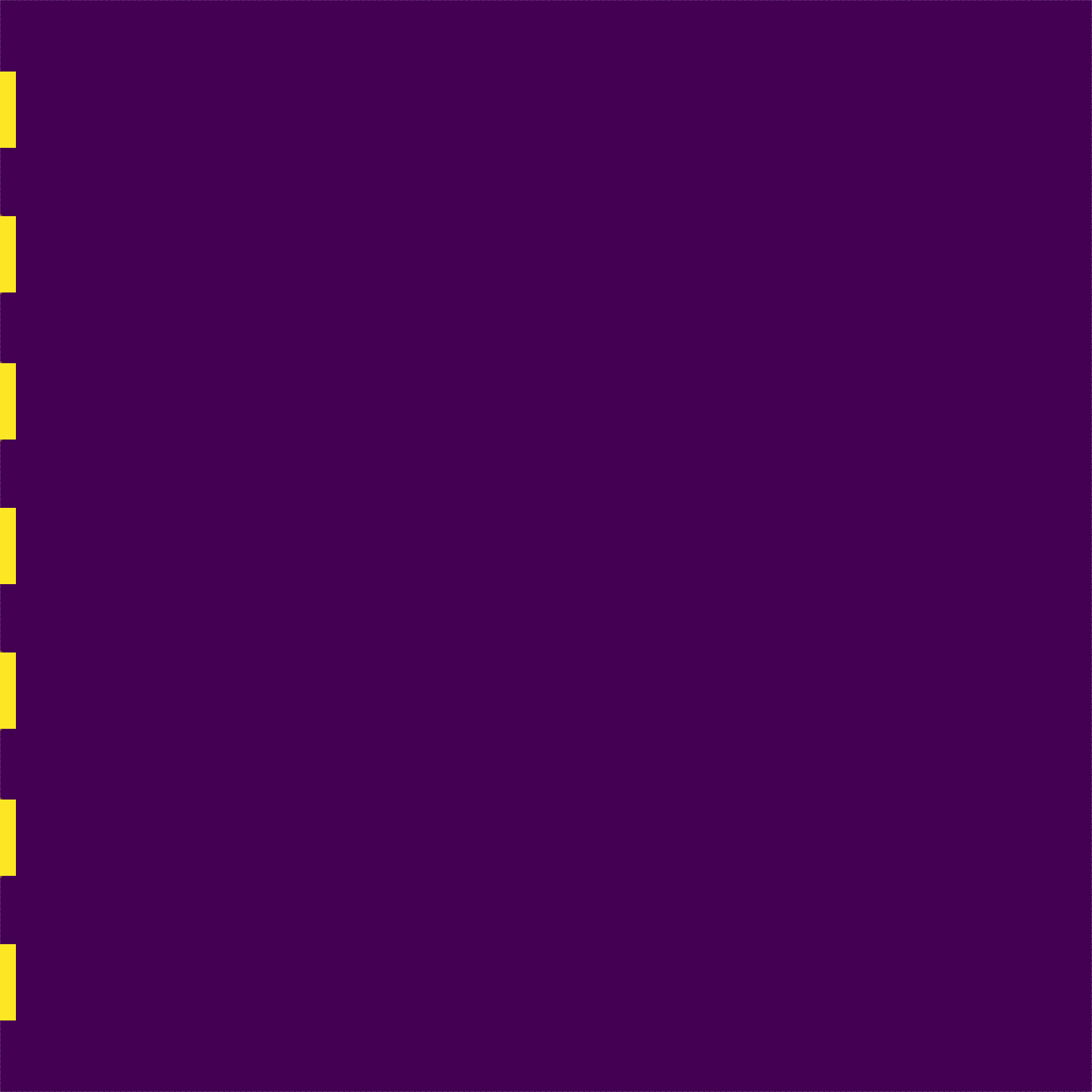};
	\end{axis}
	\node at (1.23,-0.35) {\scriptsize (c)};
	\draw [arrow, thick, black!60] (-0.6,0.25) to (0.15,0.25);
	\draw [arrow, thick, black!60] (-0.6,0.575) to (0.15,0.575);
	\draw [arrow, thick, black!60] (-0.6,0.9) to (0.15,0.9);
	\draw [arrow, thick, black!60] (-0.6,1.225) to (0.15,1.225);
	\draw [arrow, thick, black!60] (-0.6,1.55) to (0.15,1.55);
	\draw [arrow, thick, black!60] (-0.6,1.875) to (0.15,1.875);
	\draw [arrow, thick, black!60] (-0.6,2.2) to (0.15,2.2);
	\end{tikzpicture}
	\begin{tikzpicture}
	\begin{axis}[
	hide axis,
	colorbar,
	colorbar style={font=\footnotesize,width=0.15cm,ytick = {0,1}, xshift=-0.15cm},
	colormap name = viridis,
	width=4.02cm,
	height=4.02cm,
	point meta max=1,
	point meta min=0,
	xmajorticks=false,
	xmin=0, xmax=8,
	ymajorticks=false,
	ymin=0, ymax=8,
	]
	\addplot graphics [includegraphics cmd=\pgfimage,xmin=0, xmax=8, ymin=0, ymax=8] {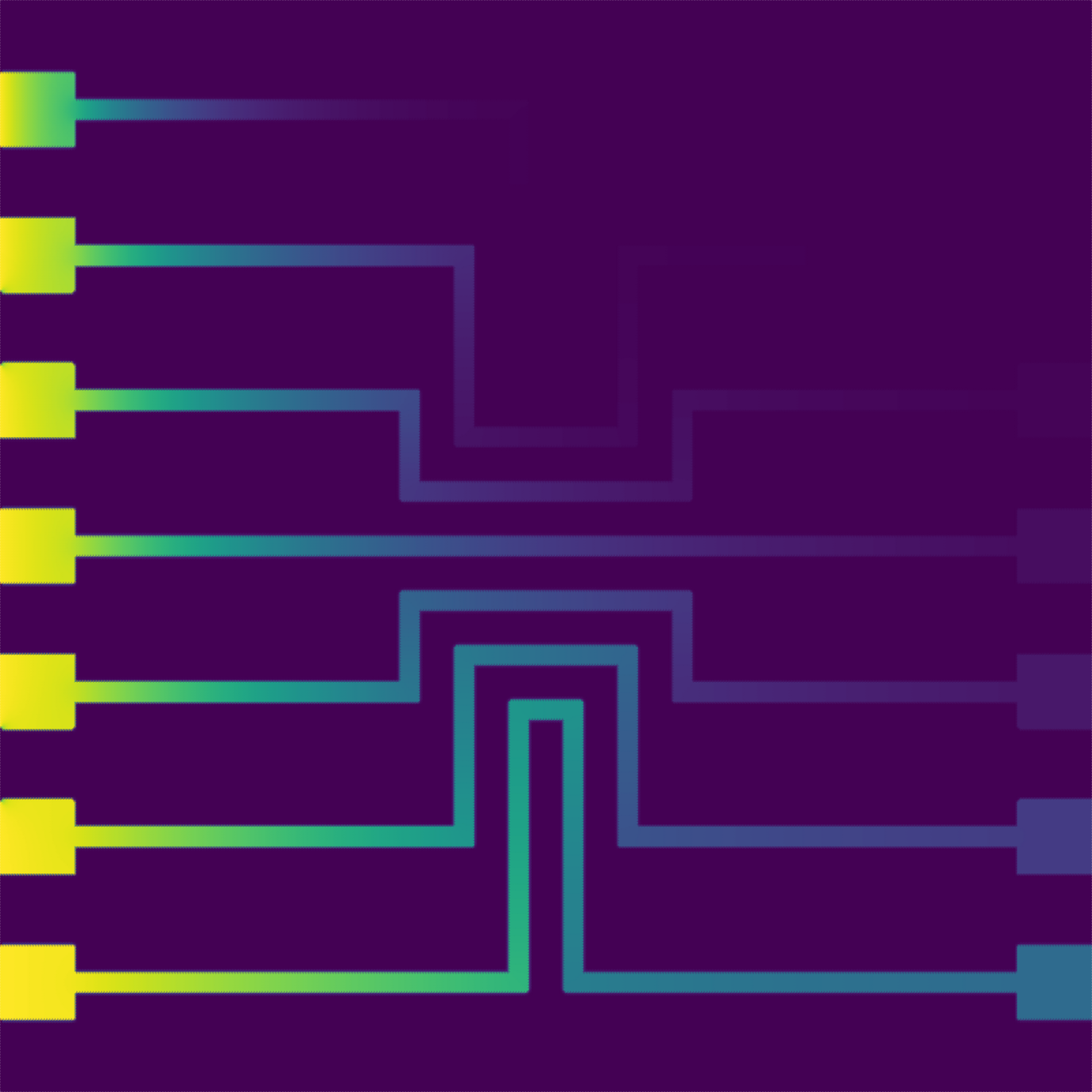};
	\end{axis}
	\node at (1.23,-0.35) {\scriptsize (d)};
	\end{tikzpicture}
	\caption{(a) Non-overlapping domain decomposition into $8\times 8$ coarse subdomains. (b) Parametric permeability $\kappa_\mu$ including channels with a conductivity ranging from $10^4$ to $10^6$. (c) Dirichlet boundary values. (d) Solution of diffusion-reaction problem for channel configuration shown in (b).} \label{fig:num_ex}
\end{figure}

We decompose the computational domain $\Omega = (0,1)^2$ into $N_H = 8 \times 8 = 64$ equally sized coarse subdomains as depicted in \cref{fig:num_ex}(a) and further discretize each subdomain with a regular quadrilateral mesh with mesh size $2^{-8}$ in both directions. We consider a parametrized heterogeneous diffusion-reaction problem of the form \cref{eq:PDE_weak} with a permeability $\kappa_\mu$ that is characterized by parametric high conductivity channels as illustrated in \cref{fig:num_ex}(b). By varying the parameter $\mu$, the conductivity within each of the seven channels is set individually in a range from $10^4$ to~$10^6$. Moreover, we choose the reaction rate $r_\mu \equiv 10^{6}$ and the source term $f_\mu \equiv 0$. We prescribe homogeneous Neumann boundary conditions on the right boundary of~$\Omega$ and Dirichlet boundary conditions on the left, top, and bottom boundary of $\Omega$ as depicted in \cref{fig:num_ex}(c).\footnote{We treat the inhomogeneous Dirichlet boundary conditions by introducing a shift function that enters the right-hand side of the problem and then solve for homogeneous Dirichlet boundary conditions on the respective part of the boundary.} In \cref{fig:num_ex}(d) the solution of the full order model (cf. problem \cref{eq:PDE_weak_DG}) for the channel configuration in \cref{fig:num_ex}(b) is visualized, showing that the distribution of mass in the channels increases with increasing conductivity. 

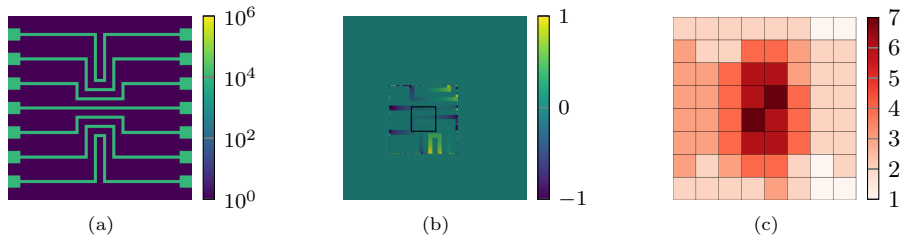
\begin{figure}[h]
	\centering
	\hfill
	\begin{tikzpicture}
	\begin{axis}[
	hide axis,
	colorbar,
	colorbar style={font=\footnotesize,width=0.15cm, ytick={0,2,4,6}, yticklabels={$10^0$,$10^2$,$10^4$,$10^6$}, xshift=-0.15cm},
	colormap name = viridis,
	width=4.02cm,
	height=4.02cm,
	point meta max=6,
	point meta min=0,
	xmajorticks=false,
	xmin=0, xmax=8,
	ymajorticks=false,
	ymin=0, ymax=8,
	]
	\addplot graphics [includegraphics cmd=\pgfimage,xmin=0, xmax=8, ymin=0, ymax=8] {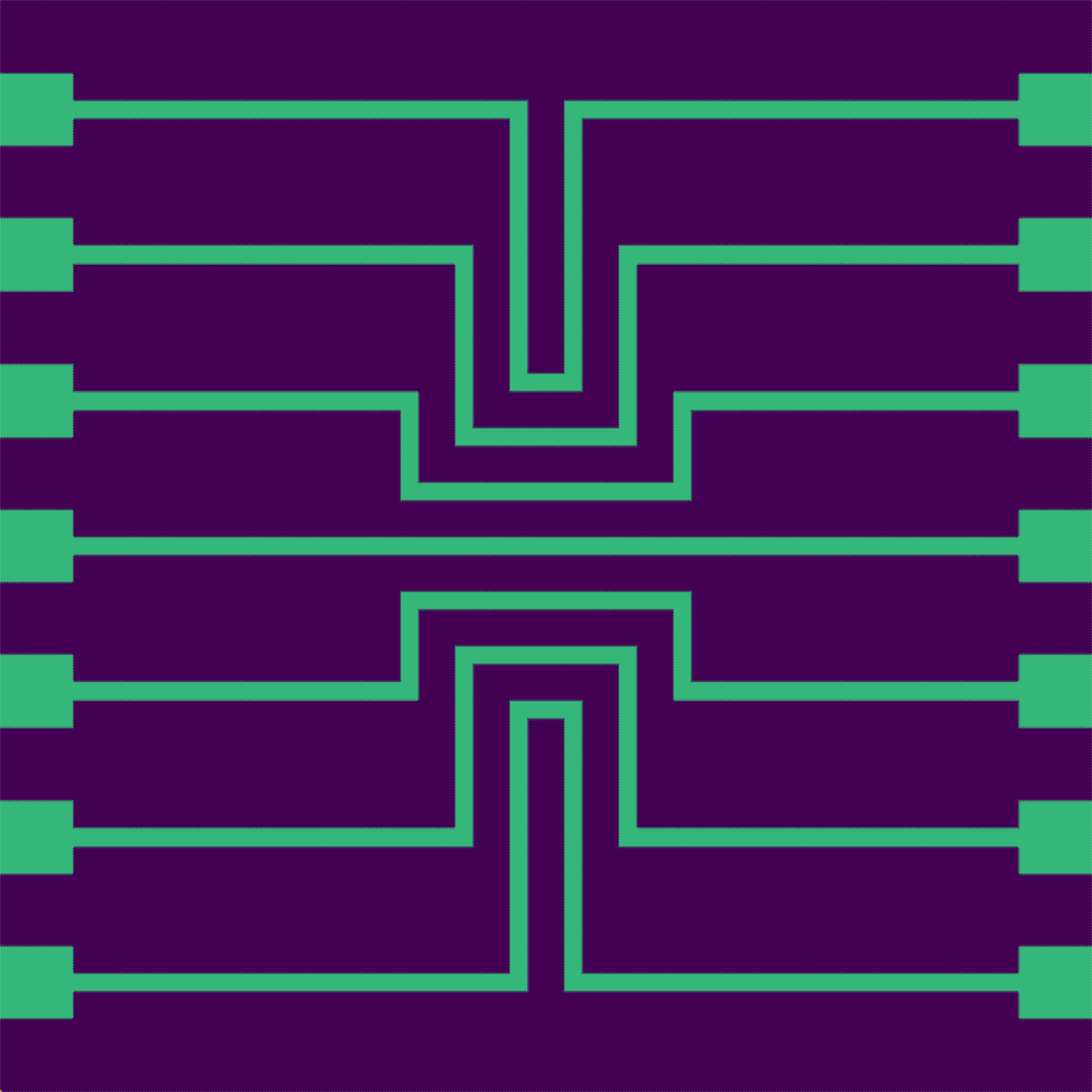};
	\end{axis}
	\node at (1.23,-0.35) {\scriptsize (a)};
	\end{tikzpicture}
	\hfill
	\begin{tikzpicture}
	\begin{axis}[
	hide axis,
	colorbar,
	colorbar style={font=\footnotesize,width=0.15cm,ytick = {-1,0,1}, xshift=-0.15cm},
	colormap name = viridis,
	width=4.02cm,
	height=4.02cm,
	point meta max=1,
	point meta min=-1,
	xmajorticks=false,
	xmin=0, xmax=8,
	ymajorticks=false,
	ymin=0, ymax=8,
	]
	\addplot graphics [includegraphics cmd=\pgfimage,xmin=0, xmax=8, ymin=0, ymax=8] {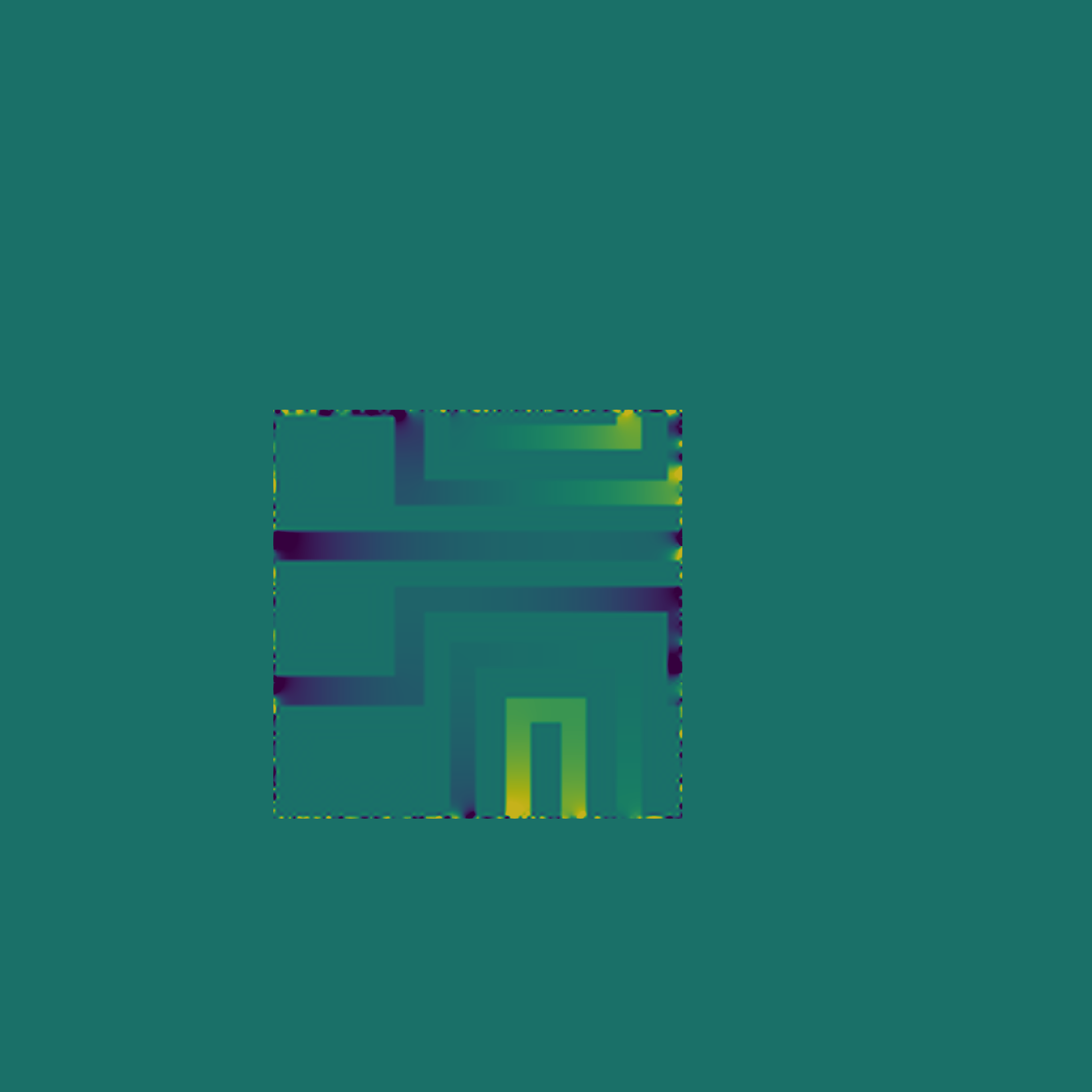};
	\end{axis}
	\draw [black] (0.905,0.905) rectangle (1.23,1.23);
	\node at (1.23,-0.35) {\scriptsize (b)};
	\end{tikzpicture}
	\hfill
	\begin{tikzpicture}
	\begin{axis}[
	colorbar,
	colorbar style={font=\small,width=0.15cm,ytick={1,2,...,7},xshift=-0.15cm},
	colormap={mymap}{[1pt]
		rgb(0pt)=(1,0.96078431372549,0.941176470588235);
		rgb(1pt)=(0.996078431372549,0.87843137254902,0.823529411764706);
		rgb(2pt)=(0.988235294117647,0.733333333333333,0.631372549019608);
		rgb(3pt)=(0.988235294117647,0.572549019607843,0.447058823529412);
		rgb(4pt)=(0.984313725490196,0.415686274509804,0.290196078431373);
		rgb(5pt)=(0.937254901960784,0.231372549019608,0.172549019607843);
		rgb(6pt)=(0.796078431372549,0.0941176470588235,0.113725490196078);
		rgb(7pt)=(0.647058823529412,0.0588235294117647,0.0823529411764706);
		rgb(8pt)=(0.403921568627451,0,0.0509803921568627)
	},
	width=4cm,
	height=4cm,
	point meta max=7,
	point meta min=1,
	xmajorticks=false,
	xmin=0, xmax=8,
	ymajorticks=false,
	ymin=0, ymax=8,
	]
	\addplot graphics [includegraphics cmd=\pgfimage,xmin=0, xmax=8, ymin=0, ymax=8] {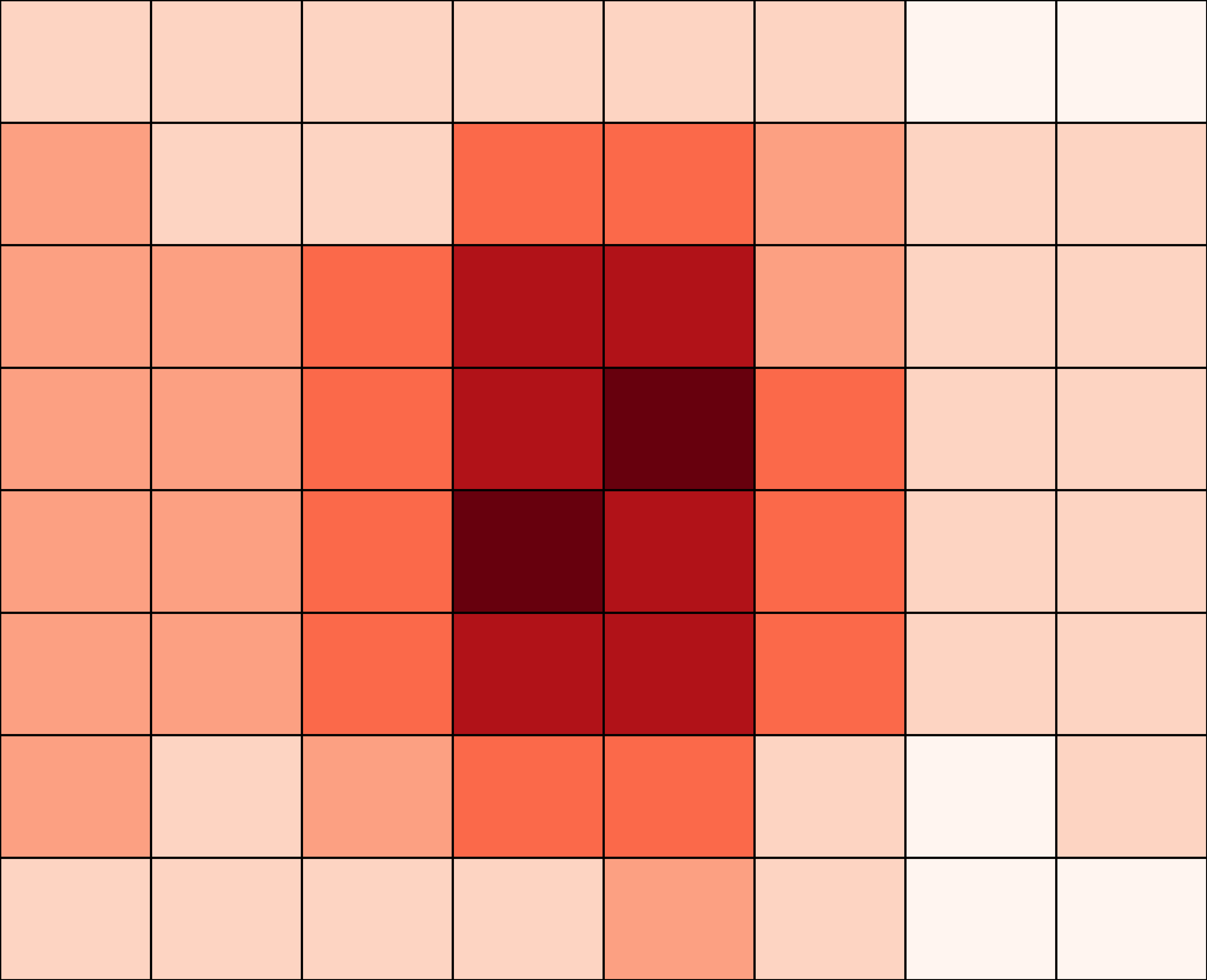};
	\end{axis}
	\node at (1.23,-0.35) {\scriptsize (c)};
	\end{tikzpicture}
	\hfill
	\caption{(a) Permeability $\kappa_{\mu_\text{train}}$ for training parameter $\mu_\text{train}$. (b) Local solution with random boundary values. (c) Number of offline computed local basis functions in randomized range finder algorithm \cite{BuhSme18} with accuracy $\texttt{tol}=10^{-2}$ and failure probability $\varepsilon_{\tt fail}=10^{-15}$.} \label{fig:num_ex_offline}
\end{figure}

\vspace{0.5em}
\textbf{Local offline training.} For the offline training phase, we choose one training parameter $\mu_\text{train}$ with associated permeability $\kappa_{\mu_\text{train}}$ depicted in \cref{fig:num_ex_offline}(a). In each coarse subdomain $T\in\mathcal{T}_H$, we then run the randomized range finder algorithm from \cite{BuhSme18} outlined in \cref{subsec:offline} for $\mu_\text{train}$, accuracy $\texttt{tol} =10^{-2}$, and failure probability $\varepsilon_{\texttt{fail}} = 10^{-15}$ to initialize the local reduced solution space $V_\text{rb}(T)$. We thus obtain $\mathbb{P}(\Vert P_{O_T \rightarrow T}^{\mu_\text{train}} - \text{proj}_{V_\text{rb}(T)} P_{O_T \rightarrow T}^{\mu_\text{train}}\Vert \leq 10^{-2}) > (1 - 10^{-15})$, cf. \cref{subsec:offline}. In \cref{fig:num_ex_offline}(c) we observe that we only require at most seven local basis functions in each subdomain $T$ to satisfy the latter property. Moreover, in \cref{fig:num_ex_offline}(b) a local solution with random Dirichlet boundary values on an oversampling domain as computed within the randomized range finder algorithm and its restriction to the associated local subdomain is depicted. Finally, we highlight that the computations for all coarse subdomains are independent of each other and can thus be easily parallelized.

\begin{figure}[t]
	\flushleft
	\begin{tikzpicture}
	\begin{axis}[
	hide axis,
	width=4.02cm,
	height=4.02cm,
	point meta max=1,
	point meta min=0,
	xmajorticks=false,
	xmin=0, xmax=8,
	ymajorticks=false,
	ymin=0, ymax=8,
	]
	\addplot graphics [includegraphics cmd=\pgfimage,xmin=0, xmax=8, ymin=0, ymax=8] {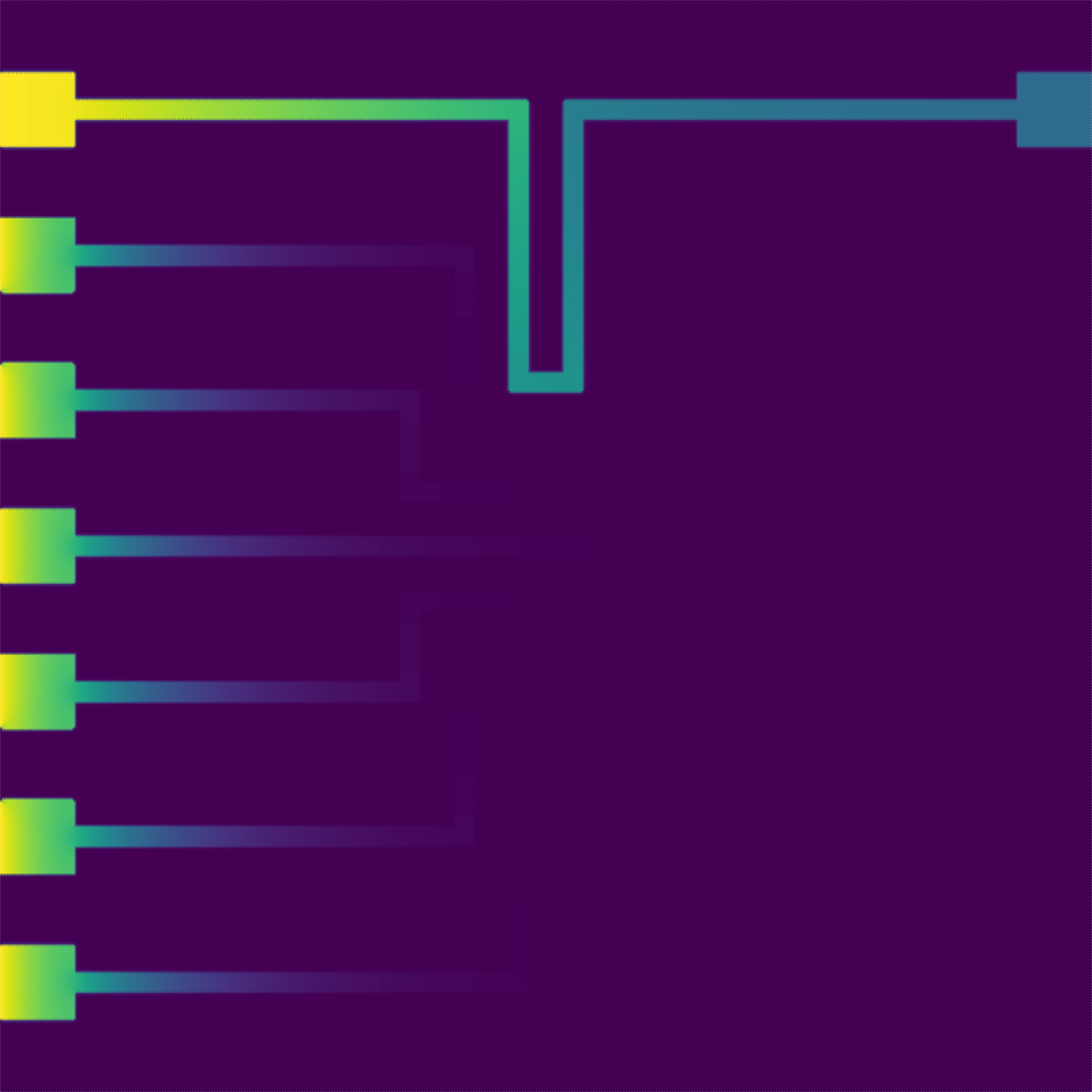};
	\end{axis}
	\end{tikzpicture}
	\hspace{0.69cm}
	\begin{tikzpicture}
	\begin{axis}[
	hide axis,
	width=4.01cm,
	height=4.01cm,
	point meta max=1,
	point meta min=0,
	xmajorticks=false,
	xmin=0, xmax=8,
	ymajorticks=false,
	ymin=0, ymax=8,
	]
	\addplot graphics [includegraphics cmd=\pgfimage,xmin=0, xmax=8, ymin=0, ymax=8] {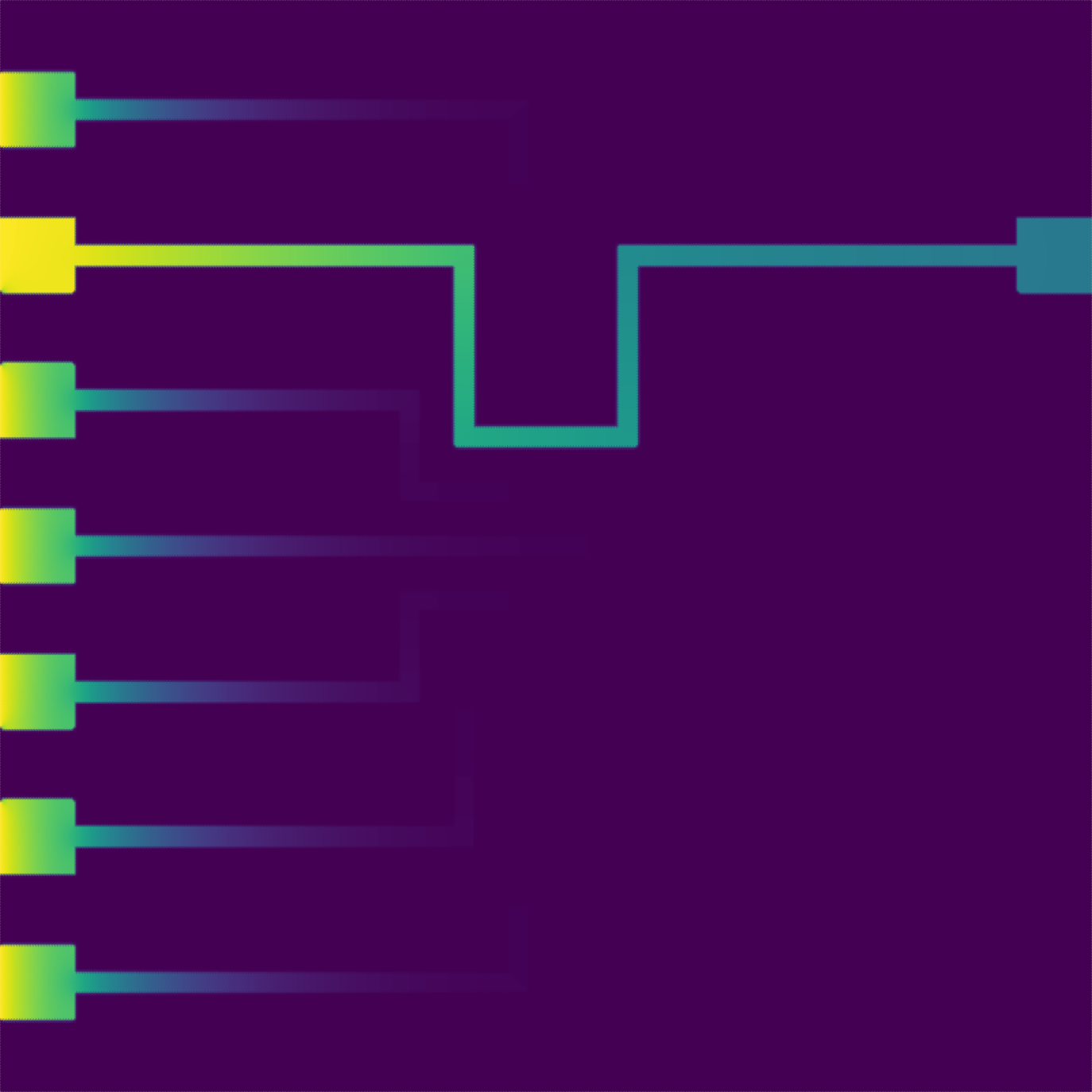};
	\end{axis}
	\end{tikzpicture}
	\hspace{0.6cm}
	\begin{tikzpicture}
	\begin{axis}[
	hide axis,
	width=4.01cm,
	height=4.01cm,
	point meta max=1,
	point meta min=0,
	xmajorticks=false,
	xmin=0, xmax=8,
	ymajorticks=false,
	ymin=0, ymax=8,
	]
	\addplot graphics [includegraphics cmd=\pgfimage,xmin=0, xmax=8, ymin=0, ymax=8] {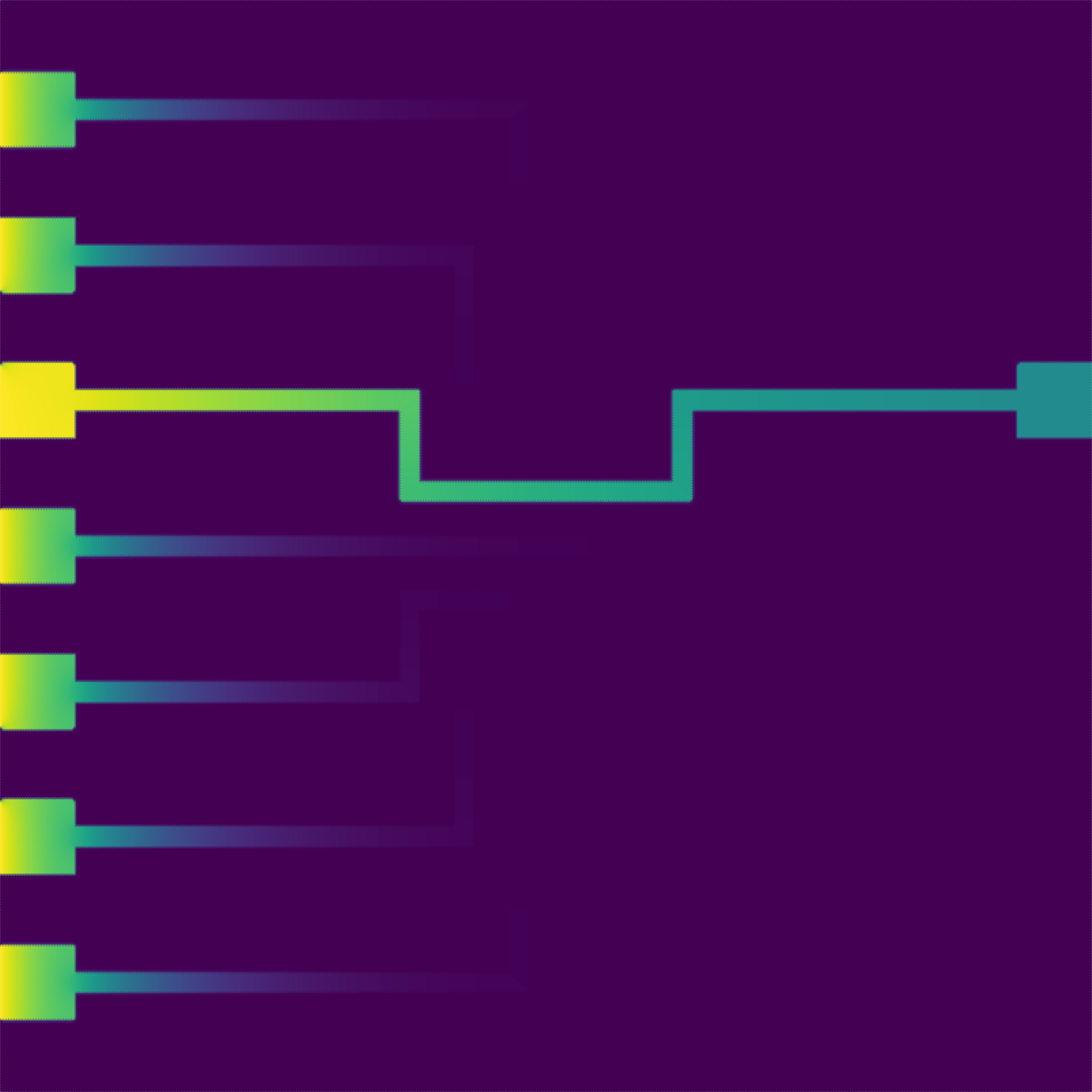};
	\end{axis}
	\end{tikzpicture}
	\hspace{0.6cm}
	\begin{tikzpicture}
	\begin{axis}[
	hide axis,
	width=4.01cm,
	height=4.01cm,
	point meta max=1,
	point meta min=0,
	xmajorticks=false,
	xmin=0, xmax=8,
	ymajorticks=false,
	ymin=0, ymax=8,
	]
	\addplot graphics [includegraphics cmd=\pgfimage,xmin=0, xmax=8, ymin=0, ymax=8] {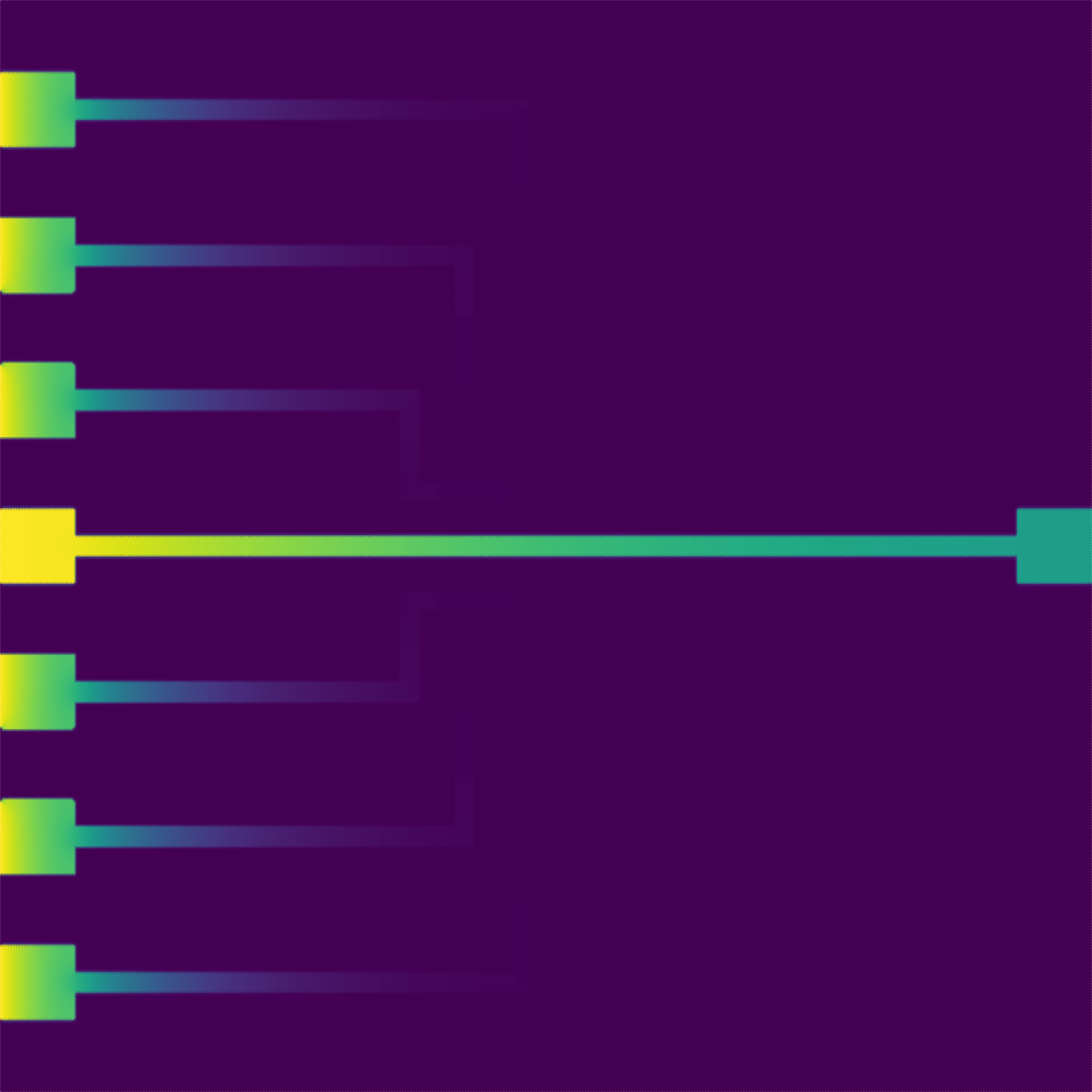};
	\end{axis}
	\end{tikzpicture}
	\vspace{-0.1cm}
	\\
	\begin{tikzpicture}
	\begin{axis}[
	colorbar,
	colorbar style={font=\small,width=0.2cm,ytick={0,1,2,3,4,5},xshift=-0.1cm},
	colormap={mymap}{[1pt]
		rgb(0pt)=(1,0.96078431372549,0.941176470588235);
		rgb(1pt)=(0.996078431372549,0.87843137254902,0.823529411764706);
		rgb(2pt)=(0.988235294117647,0.733333333333333,0.631372549019608);
		rgb(3pt)=(0.988235294117647,0.572549019607843,0.447058823529412);
		rgb(4pt)=(0.984313725490196,0.415686274509804,0.290196078431373);
		rgb(5pt)=(0.937254901960784,0.231372549019608,0.172549019607843);
		rgb(6pt)=(0.796078431372549,0.0941176470588235,0.113725490196078);
		rgb(7pt)=(0.647058823529412,0.0588235294117647,0.0823529411764706);
		rgb(8pt)=(0.403921568627451,0,0.0509803921568627)
	},
	width=4cm,
	height=4cm,
	point meta max=5,
	point meta min=0,
	xmajorticks=false,
	xmin=0, xmax=8,
	ymajorticks=false,
	ymin=0, ymax=8,
	xlabel style = {yshift=0.95cm},
	]
	\addplot graphics [includegraphics cmd=\pgfimage,xmin=0, xmax=8, ymin=0, ymax=8] {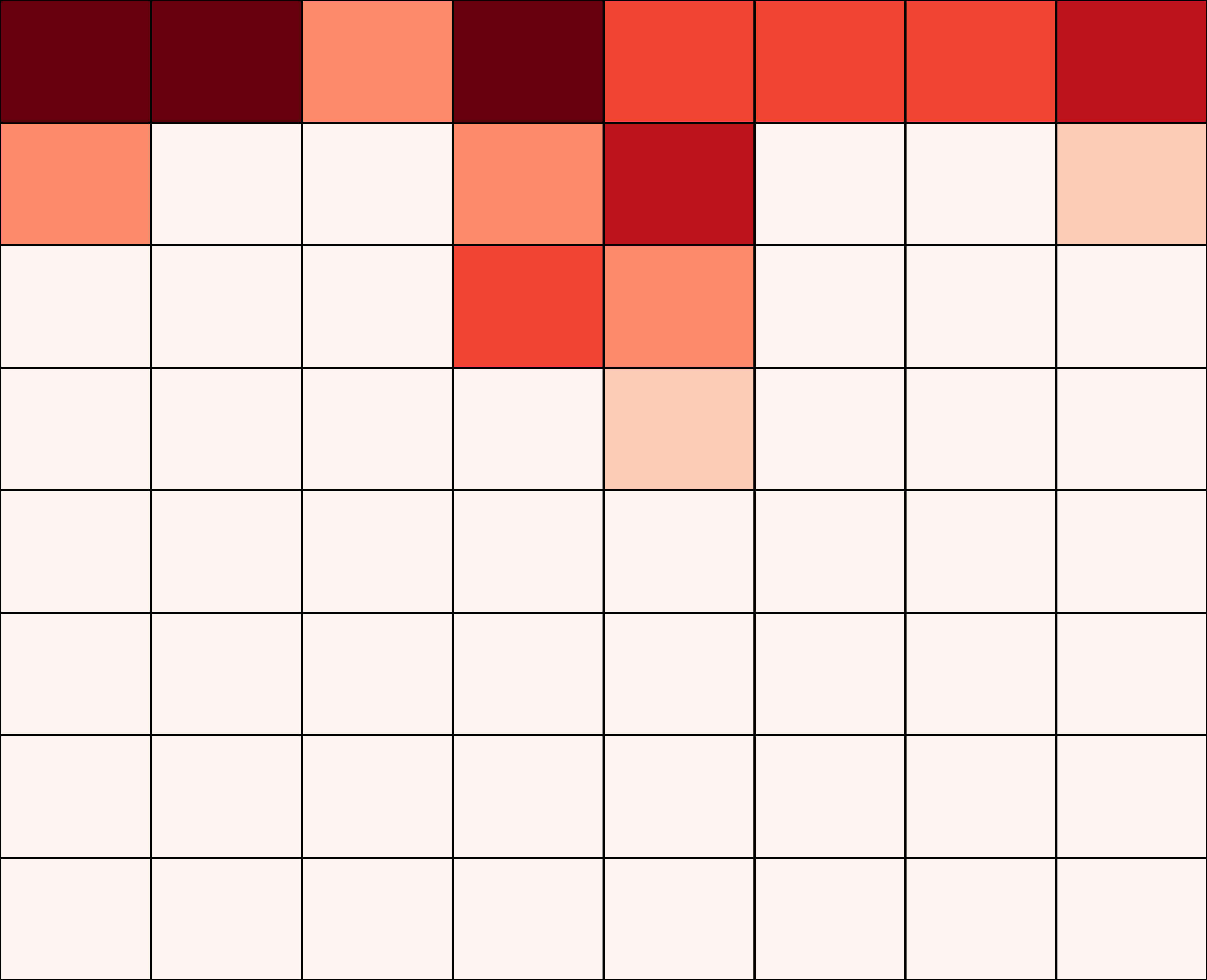};
	\end{axis}
	\end{tikzpicture}
	\hfill
	\begin{tikzpicture}
	\begin{axis}[
	colorbar,
	colorbar style={font=\small,width=0.15cm,ytick={0,1,2,3},xshift=-0.15cm},
	colormap={mymap}{[1pt]
		rgb(0pt)=(1,0.96078431372549,0.941176470588235);
		rgb(1pt)=(0.996078431372549,0.87843137254902,0.823529411764706);
		rgb(2pt)=(0.988235294117647,0.733333333333333,0.631372549019608);
		rgb(3pt)=(0.988235294117647,0.572549019607843,0.447058823529412);
		rgb(4pt)=(0.984313725490196,0.415686274509804,0.290196078431373);
		rgb(5pt)=(0.937254901960784,0.231372549019608,0.172549019607843);
		rgb(6pt)=(0.796078431372549,0.0941176470588235,0.113725490196078);
		rgb(7pt)=(0.647058823529412,0.0588235294117647,0.0823529411764706);
		rgb(8pt)=(0.403921568627451,0,0.0509803921568627)
	},
	width=4cm,
	height=4cm,
	point meta max=3,
	point meta min=0,
	xmajorticks=false,
	xmin=0, xmax=8,
	ymajorticks=false,
	ymin=0, ymax=8,
	xlabel style = {yshift=0.95cm},
	]
	\addplot graphics [includegraphics cmd=\pgfimage,xmin=0, xmax=8, ymin=0, ymax=8] {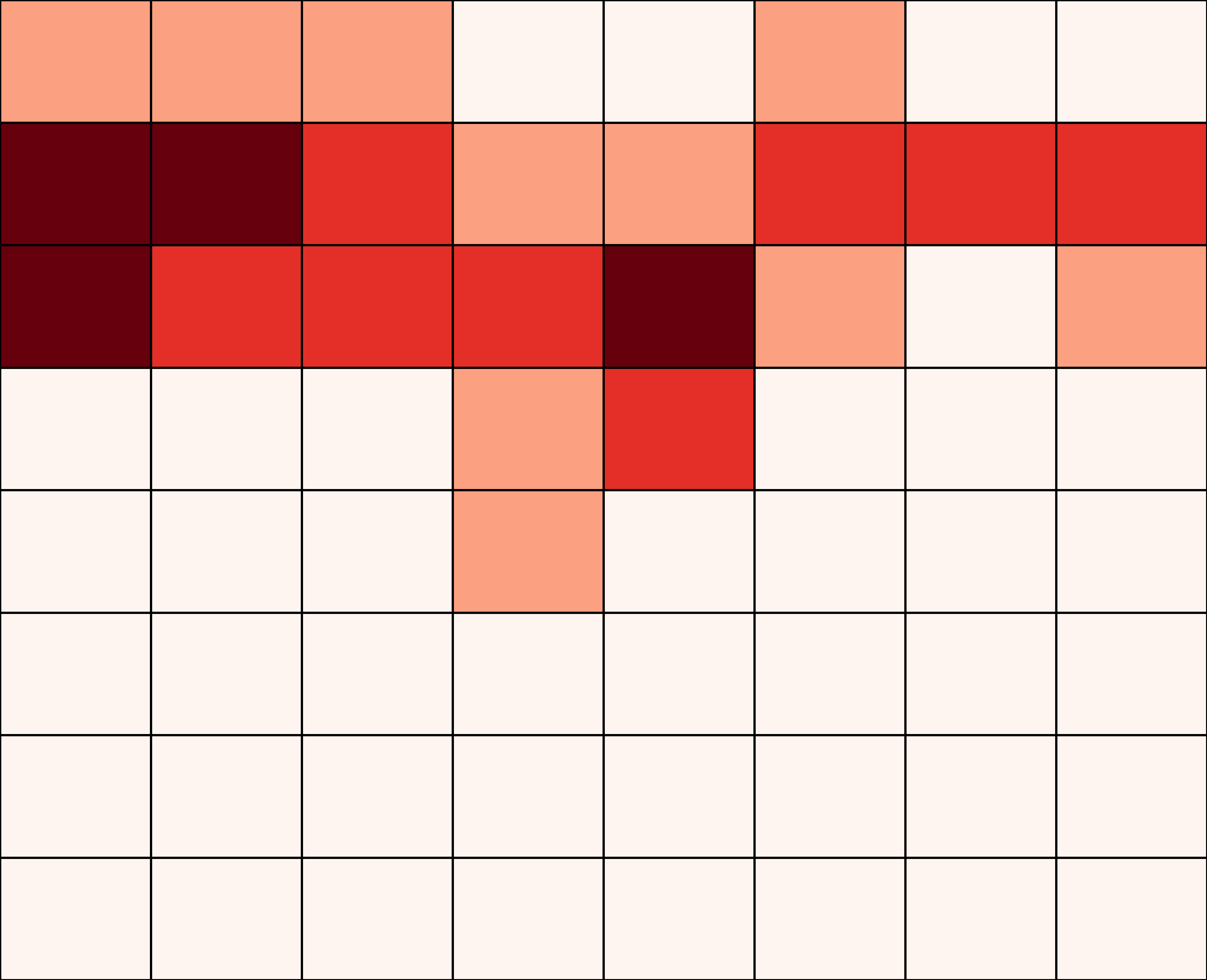};
	\end{axis}
	\end{tikzpicture}
	\hfill
	\begin{tikzpicture}
	\begin{axis}[
	colorbar,
	colorbar style={font=\small,width=0.15cm,ytick={0,1,2,3,4},xshift=-0.15cm},
	colormap={mymap}{[1pt]
		rgb(0pt)=(1,0.96078431372549,0.941176470588235);
		rgb(1pt)=(0.996078431372549,0.87843137254902,0.823529411764706);
		rgb(2pt)=(0.988235294117647,0.733333333333333,0.631372549019608);
		rgb(3pt)=(0.988235294117647,0.572549019607843,0.447058823529412);
		rgb(4pt)=(0.984313725490196,0.415686274509804,0.290196078431373);
		rgb(5pt)=(0.937254901960784,0.231372549019608,0.172549019607843);
		rgb(6pt)=(0.796078431372549,0.0941176470588235,0.113725490196078);
		rgb(7pt)=(0.647058823529412,0.0588235294117647,0.0823529411764706);
		rgb(8pt)=(0.403921568627451,0,0.0509803921568627)
	},
	width=4cm,
	height=4cm,
	point meta max=4,
	point meta min=0,
	xmajorticks=false,
	xmin=0, xmax=8,
	ymajorticks=false,
	ymin=0, ymax=8,
	xlabel style = {yshift=0.95cm},
	]
	\addplot graphics [includegraphics cmd=\pgfimage,xmin=0, xmax=8, ymin=0, ymax=8] {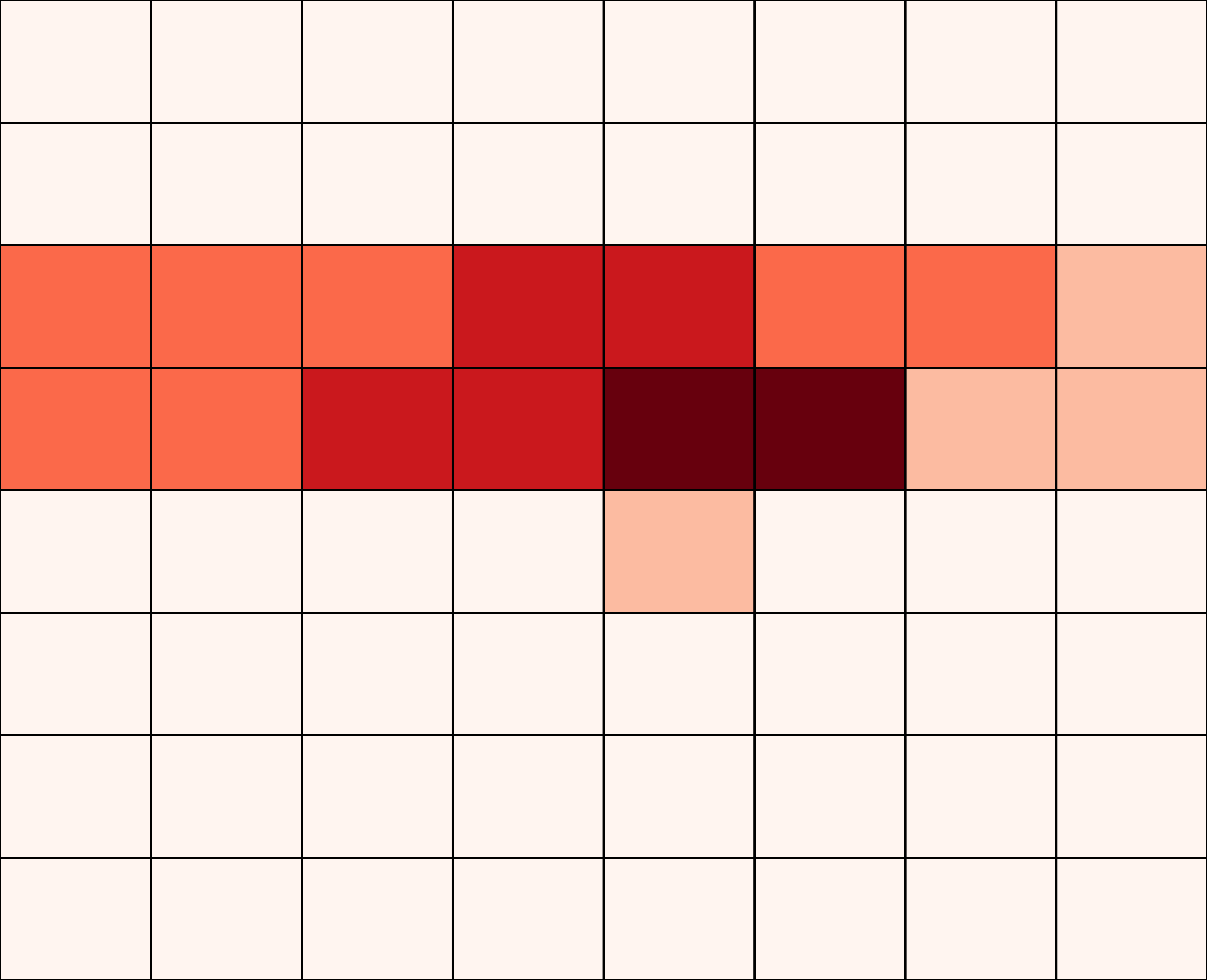};
	\end{axis}
	\end{tikzpicture}
	\hfill
	\begin{tikzpicture}
	\begin{axis}[
	colorbar,
	colorbar style={font=\small,width=0.15cm,ytick={0,1,2,3},xshift=-0.15cm},
	colormap={mymap}{[1pt]
		rgb(0pt)=(1,0.96078431372549,0.941176470588235);
		rgb(1pt)=(0.996078431372549,0.87843137254902,0.823529411764706);
		rgb(2pt)=(0.988235294117647,0.733333333333333,0.631372549019608);
		rgb(3pt)=(0.988235294117647,0.572549019607843,0.447058823529412);
		rgb(4pt)=(0.984313725490196,0.415686274509804,0.290196078431373);
		rgb(5pt)=(0.937254901960784,0.231372549019608,0.172549019607843);
		rgb(6pt)=(0.796078431372549,0.0941176470588235,0.113725490196078);
		rgb(7pt)=(0.647058823529412,0.0588235294117647,0.0823529411764706);
		rgb(8pt)=(0.403921568627451,0,0.0509803921568627)
	},
	width=4cm,
	height=4cm,
	point meta max=3,
	point meta min=0,
	xmajorticks=false,
	xmin=0, xmax=8,
	ymajorticks=false,
	ymin=0, ymax=8,
	xlabel style = {yshift=0.95cm},
	]
	\addplot graphics [includegraphics cmd=\pgfimage,xmin=0, xmax=8, ymin=0, ymax=8] {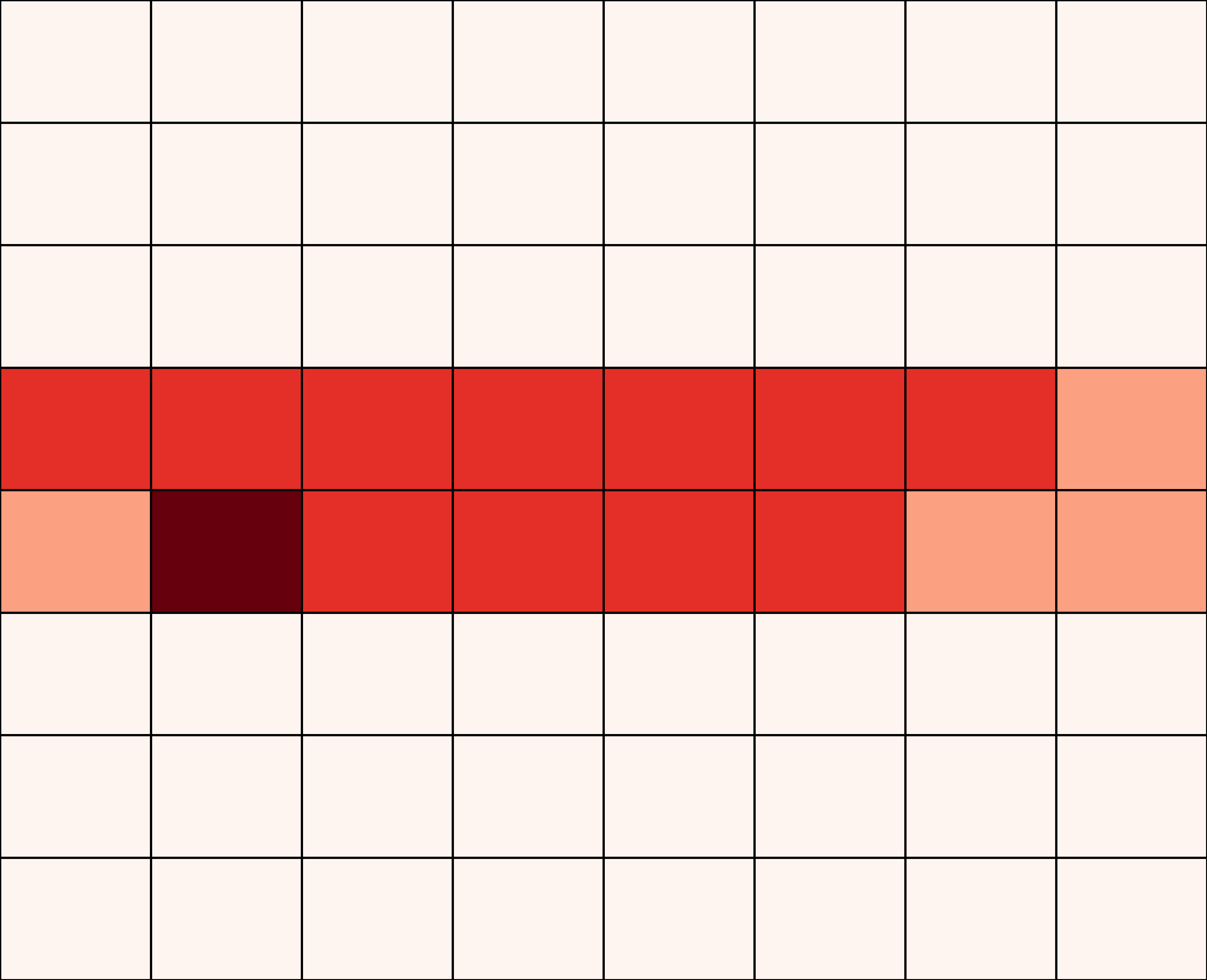};
	\end{axis}
	\end{tikzpicture}
	\caption{Top: (Full order) solutions for online requested parameters. Bottom: Number of online computed local basis functions based on localized error indicator such that relative energy error is less than $10^{-3}$.}	\label{fig:num_ex_online}
\end{figure}

\textbf{Local online enrichment with localized error control.} In the online phase, we use the localized error indicator $\delta_{T,\text{loc}}$ derived in \cref{subsec:offline} to locally investigate the accuracy of the reduced solution for new requested parameters of interest and enrich the reduced solution space locally where necessary. The four parameters chosen online in this experiment vary the permeability in such a way that the permeability in one of the channels is much higher than in the other channels, see \cref{fig:num_ex_online} (top) for a visualization of the associated (full order) solutions.

In each of the four cases, we choose the offline computed basis (cf. \cref{fig:num_ex_offline}(c) and the preceding paragraph) as the initial reduced basis. For the online enrichment procedure, we then adaptively select subdomains $T\in\mathcal{T}_H$ with the largest local error indicator $\delta_{T,\text{loc}}$ until the sum of local error indicators (of selected subdomains) exceeds $50\%$ of the global error estimate. In other words, we choose local subdomains with the worst estimated local error that in sum contribute to at least $50\%$ of the estimated global error. We then enrich the local reduced approximation space $V_\text{rb}(T)$ for all selected subdomains $T$ as outlined in \cref{subsec:offline}. In this way, we iteratively proceed until the relative energy error between the full and reduced order solution\footnote{In practice, one would employ the global error estimator \cref{eq:error_estimator} within the stopping criterion. We here consider the actual error between full and reduced order solution to demonstrate the potential of the proposed approach.} is less than $10^{-3}$.

In \cref{fig:num_ex_online} (bottom) we observe in each of the four cases that the location of the online computed local basis functions exactly aligns with the location of the channel whose permeability is changed due to the online requested parameter. We thus infer that in this numerical experiment the localized error estimator accurately detects local changes. Moreover, we emphasize that both the error estimation and the solution of the online enrichment problems require only local computations that can in addition easily be parallelized.

\section{Conclusions}
\label{sec:conclusion}

In this contribution, we proposed an adaptive localized model order reduction framework to efficiently approximate solutions of parametrized multi and/or large scale problems. The approach is based on both localized training and adaptive local enrichment. Certification and adaptivity is achieved by exploiting localized residual-type a posteriori error estimates, where localization is obtained using an abstract localization principle for dual norms.

We conjecture that the approach offers great potential to be employed, for instance, within optimization or inverse problems as it allows to significantly reduce computational costs and exploit parallelism on modern computer architectures.

\end{document}